\documentclass[12pt,intlim,righttag]{article}
\usepackage{amsmath,amsthm}
\usepackage{amssymb}
\usepackage{amsfonts}

\newcommand{\esssup}{\operatornamewithlimits{ess\,sup}}

\newcommand{\la}{\langle}
\newcommand{\ra}{\rangle}
\newcommand{\wt}{\widetilde}

\newcommand{\ld}{\lambda}

\newcommand{\td}{\tilde}

\newcommand{\be}{\begin}
\newcommand{\ee}{\end}
\newcommand{\lbl}{\label}

\newcommand\beq{\begin{equation}}
\newcommand\eeq{\end{equation}}
\newcommand{\beaa}{\begin{eqnarray*}}
\newcommand{\eeaa}{\end{eqnarray*}}

\theoremstyle{Theorem}

\theoremstyle{corollary}

\theoremstyle{remark}

\theoremstyle{definition}

\begin{document}
\title{Connections between a system of Forward-Backward SDEs and Backward Stochastic PDEs related to the utility maximization problem
 }

\author{M. Mania and R. Tevzadze}

\date{~}
\maketitle

\begin{abstract}
{Abstract. Connections between a system of Forward-Backward SDEs and Backward Stochastic PDEs related to the utility maximization
 problem is established. Besides, we derive another version of FBSDE of the same problem and prove an existence of a solution}

\end{abstract}

\bigskip


\


\

\section{Introduction}

\

We consider a financial market model, where the dynamics of asset prices is described by the continuous $R^d$-valued continuous semimartingale $S$  defined on a
complete  probability space $(\Omega,\mathcal{F}, P)$  with filtration $F=({F}_t,t\in[0,T])$ satisfying the usual conditions, where $\mathcal{F}=F_T$ and $T<\infty$.
We work with discounted terms, i.e. the bond is assumed to be a constant.

Let $U=U(x):R\to R$ be a utility function taking finite values at all points of real line $R$ such that $U$ is continuously
differentiable, increasing, strictly concave and satisfies the Inada conditions
\begin{equation}\label{inada}
    U'(\infty)=\lim_{x\to\infty}U'(x)=0,\quad U'(-\infty)=\lim_{x\to -\infty}U'(x)=\infty.
\end{equation}
We also assume that $U$ satisfies the condition of reasonable  asymptotic elasticity (see \cite{KSch} and \cite{S3} for a detailed discussion of these conditions), i.e.
\begin{equation}\label{ae}
    \limsup_{x\to\infty}\frac{xU'(x)}{U(x)}<1,\quad \liminf_{x\to -\infty}\frac{xU'(x)}{U(x)}>1.
\end{equation}

 For the utility function $U$ we denote by ${\widetilde U}$ its convex conjugate
\begin{equation}
    \widetilde U(y)=\sup_x(U(x)-xy),\quad  y>0.
\end{equation}

Denote by ${\mathcal M}^e$ (resp. ${\mathcal M}^a$)
the set of probability  measures $Q$ equivalent (resp. absolutely continuous) with respect to $P$
such that $S$ is a local martingale under $Q$.

Let  ${\cal M}^a_U$ (resp. ${\cal M}^e_U$) be the convex set of probability measures
$Q\in{\cal M}^a $ (resp. ${\cal M}^e$) such that
\begin{equation}\label{mq}
E{\widetilde U}\big(\frac{dQ_T}{dP_T}\big)<\infty.
\end{equation}
It follows from proposition 4.1 of \cite{S1} that ( \ref{mq}) implies $E{\widetilde U}\big(y\frac{dQ_T}{dP_T}\big)<\infty$ for any $y>0$.

Throughout the paper we assume that
\begin{equation}\label{mm}
{\mathcal M}^e_U\neq\emptyset.
\end{equation}

The wealth process, determined by a self-financing
trading strategy $\pi$ and initial capital $x$,  is defined as a stochastic integral
$$
X^{x,\pi}_t=x+\int_0^t\pi_udS_u,\;\;\;0\le t\le T.
$$

We consider the utility maximization problem with random endowment $H$, where $H$ is a liability that the agent must deliver at terminal time $T$.
$H$ is an $F_T$-measurable random variable
which  for simplicity is assumed to be bounded (one can use also weaker assumption 1.6 from \cite{OZ}).
The value function $V(x)$ associated to the problem is defined by
\begin{equation}\label{ut1}
    V(x)=\sup_{\pi\in\Pi_x} E\bigg[ U\bigg(x+\int_0^T\pi_u\,dS_u+H \bigg)\bigg],
 \end{equation}
where $\Pi_x$ is a class of strategies which  (following \cite{S3} and \cite{OZ}) we define
as  the class  of predictable $S$- integrable processes
$\pi$ such that $U(x+(\pi\cdot S)_T+H)\in L^1(P)$ and $\pi\cdot S$ is
a supermartingale
under each $Q\in {\cal M}^a_U$.

The dual problem to (\ref{ut1}) is
\begin{equation}\label{dualv}
\wt V(y)=\inf_{Q\in{\cal M}^e_U}E[\wt U(y\rho_T^Q)+y\rho_T^Q],\;\;\;y>0,
\end{equation}
where $\rho^Q_t=dQ_t/dP_t$ is the density process of the measure $Q\in{\cal M}^e$ relative to the basic measure $P$.

It was shown in \cite{OZ} that under assumptions (2) and (5) an optimal strategy $\pi(x)$ in the class $\Pi_x$ exists. There exists also
an optimal martingale measure $Q(y)$ to the problem (\ref{dualv}), called the minimax martingale measure and by $\rho^*=(\rho_t(y), t\in [0,T])$ we
denote the density process of this measure relative to the measure $P$.

It follows also  from \cite{OZ} that under  assumptions
(2) and (5) optimal solutions $\pi(x)\in\Pi$ and $Q(y)\in{\cal
M}^e_U$   are
related as
\begin{equation}\label{Z}
U^{'}\left(x+\int_0^T\pi_u(x)dS_u+H\right)=y \rho_T(y),\;\;\; P-a.s.
\end{equation}

The continuity of $S$ and the existence of an equivalent martingale measure imply that the structure condition is satisfied, i.e. $S$ admits the decomposition
$$  S_t=M_t+\int_0^td\langle M\rangle_s\lambda_s,\quad  \int_0^t \lambda^T_s\,d\langle M\rangle_s\lambda_s <\infty      $$
for all $t$ $P$-a.s., where $M$ is a continuous local martingale and $\lambda$ is a predictable process.  The sign ${}^T$ here denotes the transposition.

Let us introduce a dynamic value function of the problem (\ref{ut1}) defined as
\begin{equation}\label{ut5}
    V(t,x)=\esssup_{\pi\in\Pi_x} E\bigg(U\bigg(x+\int_t^T\pi_u\,dS_u +H\bigg)\;\Big| \; {F}_t\bigg).
\end{equation}

It is well known that for any $x\in R$ the process
$(V(t,x),t\in [0,T])$ is a supermartingale admitting an RCLL (right-continuous with left limits) modification.

Therefore, using the Galchouk--Kunita--Watanabe (GKW) decomposition, the value function is represented as
$$  V(t,x)=V(0, x)- A(t,x)+\int_0^t\psi(s,x)\,dM_s+L(t,x),        $$
where for any $x\in R$ the process $A(t,x)$ is  increasing  and $L(t,x)$ is a local martingale orthogonal to $M$.

{\bf Definition 1.} We shall say that $(V(t,x),t\in[0,T])$ is a regular family of semimartingales if

a)  $V(t,x)$ is  two-times continuously differentiable at $x$  $P$- a.s. for any $t\in [0,T]$,

b) for any $x\in R$  the process $V(t,x)$ is a special semimartingale with bounded variation part absolutely continuous with respect to
an increasing predictable process  $(K_t, t\in[0,T])$, i.e.
$$  A(t,x)=\int_0^ta(s,x)\,dK_s,       $$
for some real-valued function $a(s,x)$ which is predictable  and $K$-integrable for any $x\in R$,

c) for any $x\in R$ the process $V'(t,x)$ is a special semimartingale with the decomposition
$$  V'(t,x)=V'(0, x)- \int_0^ta'(s,x)\,dK_s+\int_0^t\psi'(s,x)\,dM_s+L'(t,x).     $$
where $a', \varphi'$ and $L'$ are partial derivatives of $a, \varphi$ and $L$ respectively.

If $F(t,x)$ is a family of semimartigales   then $\int_0^TF(ds,\xi_s)$ denotes a generalized stochastic integral, or a stochastic line integral (see \cite{K},
or \cite{Ch}).
If $F(t,x)=xG_t$, where $G_t$ is a semimartingale then  the stochastic line integral  coincides with  the usual stochastic integral
denoted by $\int_0^T\xi_sdG_s$ or $(\xi\cdot G)_T$.

It was shown in \cite{MT3, MT8,MT10} (see, e.g., Theorem 3.1 from  \cite{MT10}) that if the value function  satisfies conditions a)-c) then it solves
the following BSPDE
$$
V(t,x)= V(0,x)
$$
$$+\frac{1}{2}\int_0^t \frac{1} {V''(s,x)}(\varphi'(s,x)+\lambda(s)V'(s,x))^T\,d\langle M\rangle_s(\varphi'(s,x)+\lambda(s)V'(s,x))
$$
\begin{equation}\label{pde}
   +\int_0^t\varphi(s,x)\,dM_s+L(t,x),\quad  V(T,x)=U(x)
\end{equation}
and optimal wealth satisfies the  SDE
\begin{equation}\label{pdex}
X_t(x)=x- \int_0^t\frac{\varphi'(s,X_s(x))+
\lambda(s)V'(s,X_s(x))} {V''(s,X_s(x))} dS_s.
\end{equation}
Note that  the BSPDE (\ref{pde}), (\ref{pdex}) is of the same form  for utility functions defined on half real  line and also for random utility functions $U(\omega, x)$.

In the paper \cite{hor} a new approach was developed, where a characterization of optimal strategies to the problem (\ref{ut1}) in terms of a  system of Forward-Backward Stochastic Differential Equations
(FBSDE) in the Brownian framework was given. The key observation was an existence of a stochastic process $Y$ with $Y_T=H$ such that $U'(X_t+Y_t)$ is a martingale.
The same approach was used in \cite{San}, where these  results were generalized  in
semimartingale setting with continuous filtration rejecting  also some technical  conditions imposed in \cite{hor}. The FBSDE for the pair
 $(X,Y)$ (where $X$ is the optimal wealth and  $Y$ the process mentioned above)
is of the form
\begin{equation}\label{fb}
Y_t=Y_0+ \int_0^{t}\big[ \lambda_s^T\frac{U'(X_s+Y_s)}{U''(X_s+Y_s)}-
\frac{1}{2}\lambda_s^T\frac{U'''(X_s+Y_s)U'(X_s+Y_s)^2}{U''(X_s+Y_s)^3}
\end{equation}
$$
+Z_s^T\big]d\langle M\rangle_s\lambda_s-\frac{1}{2}\int_0^t \frac{U'''(X_s+Y_s)}{U''(X_s+Y_s)}d\langle N\rangle_s +\int_0^tZ_sdM_s+N_t,\,\,\; Y_T=H.
$$
\begin{equation}\label{fbx}
X_t= x- \int_0^t\big(\lambda_s\frac{U'(X_s+Y_s)}{U''(X_s+Y_s)}+Z_s\big)dS_s,
 \end{equation}
where $N$ is a local martingale orthogonal to $M$.

Note that in (\cite{hor}) and (\cite{San})
an existence of a solution of FBSDE  (\ref{fb}), (\ref{fbx}) is not proved, since not all conditions of correspondsing theorems are formulated in terms of basic objects. E.g., in both papers is imposed that $E(U'(X_T^*+H))^2<\infty$ and it is not clear if an optimal strategy  satisfying this condition exists.
One our goal is to derive an other version of  FBSDE (\ref{fb}), (\ref{fbx}) and to prove an existence of a solution which will imply an existence of a solution of the system  (\ref{fb}), (\ref{fbx}).

The second   goal  is     to establish relations between  equations  BSPDE (\ref{pde}), (\ref{pdex})
and  FBSDE (\ref{fb}), (\ref{fbx}).
 Solutions of these equations
give constructions of the optimal strategy of one and the same problem, hence they should be related in some way. On the other hand  BSPDE (\ref{pdeder}),(\ref{pdex2})
can be considered as a generalization of Hamiltom-Jacobi-Bellman equation to the non Markovian case and  FBSDE (\ref{fb}), (\ref{fbx}) is linked with the stochastic maximum principle
(see \cite{hor}), although   equation (\ref{fb})- (\ref{fbx}) is not  obtained directly from the maximum principle.
It is well known that the relation between Bellman's dynamic programing and the Pontriagin's maximum principle in optimal control  is of the form
$\psi_t=V'(t,X_t)$, where  $V$ is the value function, $X$ an optimal solution and $\psi$ is an adjoint process (see, e.g. \cite{bs}, \cite{Z}).
Therefore,  somewhat similar relation between above mentioned equations should be expected. In particular,  it is shown in Theorem 1, that the first conponents of solutions of of these
equations are related by the equality
$$
Y_t=-{\wt U}'(V'(t,X_t))-X_t.
$$

In section 3 we derive other version of the
FBSDE system   (\ref{fb}), (\ref{fbx}) with decoupling field $u(t,x)=V'(t,x)-U'(x)$ (see definition 2 below), where the backward component $P_t$ is a process, such that $P_t+U'(X_t)$ is a martingale.

\section{Relations between BSPDE   (\ref{pde})-(\ref{pdex})
and  FBSDE (\ref{fb})-(\ref{fbx}) }

\

To establish relations between  equations  BSPDE (\ref{pde}), (\ref{pdex})
and  FBSDE (\ref{fb}), (\ref{fbx}) we need the following

 {\bf Definition 2 (\cite{Fr})}. The function $u(t,x)$
 is called  a  decoupling field of  the  FBSDE (\ref{fb}), (\ref{fbx})  if
\begin{equation}\label{uth}
u(T,x)=H
\end{equation}
and for any  $x\in R, s, \tau\in R_+$  such that $0\le s<\tau\le T$ the   FBSDE
\begin{equation}\label{fbs}
Y_t=u(s,x)
\end{equation}
$$
+\int_s^{t}\big( \lambda_r^T\frac{U'(X_r+Y_r)}{U''(X_r+Y_r)}-
\frac{1}{2}\lambda_r^T\frac{U'''(X_r+Y_r)U'(X_r+Y_r)^2}{U''(X_r+Y_r)^3}+Z_r^T\big)d\langle M\rangle_r\lambda_r
$$
$$
-\frac{1}{2}\int_s^t \frac{U'''(X_r+Y_r)}{U''(X_r+Y_r)}d\langle N\rangle_r +\int_s^tZ_rdM_r+ N_t-N_s,\,\,\; Y_\tau=u(\tau, X_\tau),
$$
\begin{equation}\label{fbxs}
X_t= x- \int_s^t\big(\lambda_r\frac{U'(X_r+Y_r)}{U''(X_r+Y_r)}+Z_r\big)dS_r,
 \end{equation}
has a solution
$(Y, Z, N, X)$  satisfying
\begin{equation}\label{utht}
Y_t=u(t,X_t), \;\;\;\;t\in[s,\tau].
\end{equation}
We shall say that $u(t,x)$ is a regular decoupling field if it is a regular family of semimartingales (in the sense of Definition 1).

If we differentiate equation BSPDE (\ref{pde}) at $x$ (assuming that all derivatives involved exist), we obtain the BSPDE
$$
    V'(t,x)= V'(0,x)
$$
$$
+\frac{1}{2}\int_0^t \bigg(\frac{(\varphi'(s,x)+\lambda_sV'(s,x))^T} {V''(s,x)}\,d\langle M\rangle_s(\varphi'(s,x)+\lambda_sV'(s,x))\bigg)'
$$
\begin{equation}\label{pded}
   +\int_0^t\varphi'(s,x)\,dM_s+ L'(t,x),\quad  V'(T,x)=U'(x+H).
\end{equation}
Thus, we consider the following BSPDE
$$
    V'(t,x)= V'(0,x) +\int_0^t \bigg[\frac{(V''(s, x)\lambda_s+\varphi''(s, x))^T}{V''(s, x)}
$$
$$
-\frac{1}{2}V'''(s, x)\frac{(V'(s, x)\lambda_s+\varphi'(s, x))^T}{V''(s, x)}\bigg]\,d\langle M\rangle_s(V'(s, x)\lambda_s+\varphi'(s, x))
$$
\begin{equation}\label{pdeder}
   +\int_0^t\varphi'(s,x)\,dM_s+L'(t,x),\quad  V'(T,x)=U'(x+H),
\end{equation}
where the optimal wealth satisfies the  same SDE
\begin{equation}\label{pdex2}
X_t(x)=x- \int_0^t\frac{\varphi'(s,X_s(x))+
\lambda(s)V'(s,X_s(x))} {V''(s,X_s(x))} dS_s.
\end{equation}
The FBSDE (\ref{fb}), (\ref{fbx}) is equivalent, in some sense, to BSPDE (\ref{pdeder}),(\ref{pdex2}) and
the following statement establishes a  relation between  these equations.

{\bf Theorem 1.}
Let  the utility function $U(x)$ be three-times continuously differentiable and  let  the filtration $F$  be continuous.

 a) If $V'(t,x)$ is a regular family of semimartingales and $(V'(t,x), \varphi'(t,x), L'(t,x), X_t)$  is a solution of BSPDE (\ref{pdeder}),(\ref{pdex2}), then
 the quadruple\\
$(Y_t, Z_t, N_t, X_t)$, where
\begin{equation}\label{y}
Y_t=-{\wt U}'(V'(t,X_t))-X_t,
\end{equation}
\begin{equation}\label{z}
Z_t=\lambda_t{\wt U}'(V'(t,X_t))+\frac{\varphi'(t,X_t)+\lambda_tV'(t,X_t)}{V''(t,X_t)},
\end{equation}
\begin{equation}\label{zb}
N_t=- \int_0^t{\wt U}''({V}'(s,X_s))d\big(\int_0^sL'(dr,X_r)\big),
\end{equation}
will satisfy the FBSDE (\ref{fb}), (\ref{fbx}). Moreover, the function $u(t,x)=-{\wt U}'(V'(t,x))- x$
will be the  decoupling field of  this FBSDE.

b)  Let $u(t,x)$  be a regular decoupling field  of  FBSDE (\ref{fb}), (\ref{fbx})  and let $(U'(X_t+Y_t), s\le t\le T)$ be a true martingale for every
$s\in [0,T]$. Then $(V'(t,x), \varphi'(t,x), L'(t,x), X)$ will be a solution of BSPDE (\ref{pdeder}),(\ref{pdex2}) and
following relations hold
\begin{equation}\label{y2}
V'(t,x)=U'(x+u(t,x)),\;\;\;\;\text{hence}\;\;\;\;V'(t,X_t)=U'(X_t+Y_t),
\end{equation}
\begin{equation}\label{y3}
\varphi'(t,X_t)=(Z_t+\lambda_s\frac{U'(X_t+Y_t)}{U''(X_t+Y_t)})V''(t,X_t)-\lambda_tU'(X_t+Y_t),
\end{equation}
\begin{equation}\label{y4}
\int_0^tL'(ds,X_s)=\int_0^t{U''(X_s+Y_s)}dN_s,
\end{equation}
where $ \int_0^tL'(ds,X_s)$ is a stochastic line integral with respect to the family $(L'(t, x), x\in R)$ along the process $X$.

{\it Proof.} a) It follows from BSPDE (\ref{pdeder}), (\ref{pdex2}) and from the  It\^o-Ventzel formula that $V'(t,X_t)$ is a local martingale with the decomposition
\begin{equation}\label{op}
V'(t,X_t)=V'(0,x) - \int_0^t\lambda_sV'(s,X_s)dM_s+\int_0^tL'(ds,X_s).
\end{equation}

Let $Y_t=-{\wt U}'(V'(t,X_t))-X_t$.
Since $U$ is three-times differentiable (hence so is $\tilde U$ also), $Y_t$ will be a special semimartingale and by GKW decomposition
\begin{equation}\label{dec}
Y_t=Y_0+A_t+\int_0^tZ_udM_u+N_t,
\end{equation}
where $A$ is a predictable process of finite variations and $N$ is a local martingale orthogonal to $M$.

The definition of the process $Y$,  decompositions (\ref{op}) ,  (\ref{dec}) and the It\^o formula for ${\wt U}'(V'(t,X_t))$  imply that
 \begin{equation}\label{eq}
A_t+\int_0^tZ_sdM_s+N_t=
 \end{equation}
$$
=\int_0^{t}{\wt U}''(V'(s,X_s)) V'(s,X_s)\lambda_sdM_s-\int_0^{t}{\wt U}''(V'(s,X_s)) d\big(\int_0^sL'(dr,X_r)\big)
$$
$$
-\frac{1}{2}\int_0^{t}{\wt U}'''(V'(s,X_s)) V'(s,X_s)^2\lambda^T_sd\la M\ra_s\lambda_s- \frac{1}{2}\int_0^{t}{\wt U}'''(V'(s,X_s)) d\la \int_0^.L'(dr,X_r)\ra_s
$$
$$
+\int_0^{t}\frac{\lambda_sV'(s,X_s)+\varphi'(s,X_s)}{V''(s,X_s)}dM_s+\int_0^{t}\frac{\lambda_s^TV'(s,X_s)+\varphi'(s,X_s)^T}{V''(s,X_s)}d\la M\ra_s\lambda_s,
$$
Equalizing the integrands of stochastic integrals with respect to $dM$   in (\ref{eq}) we have that $\mu^{\la M\ra>}$-a.e.
\begin{equation}\label{eq2}
Z_s=\frac{\lambda_sV'(s,X_s)+\varphi'(s,X_s)}{V''(s,X_s)}+{\wt U}''(V'(s,X_s)) V'(s,X_s)\lambda_s.
 \end{equation}
Equalizing the orthogonal martingale parts  we get $P$-a.s.
\begin{equation}\label{eq3}
N_t=-\int_0^t{\wt U}''(V'(s,X_s))d\big(\int_0^sL'(dr,X_r)\big).
 \end{equation}
Equalizing the parts of finite variations in (\ref{eq}) we have
\begin{equation}\label{eq4}
A_t=\int_0^{t}\frac{\lambda_s^TV'(s,X_s)+\varphi'(s,X_s)^T}{V''(s,X_s)}d\la M\ra_s\lambda_s
\end{equation}
$$
-\frac{1}{2}\int_0^{t}{\wt U}'''(V'(s,X_s)) V'(s,X_s)^2\lambda^T_sd\la M\ra_s\lambda_s- \frac{1}{2}\int_0^{t}{\wt U}'''(V'(s,X_s)) d\la \int_0^.L'(dr,X_r)\ra_s
$$
and  by  equalities  (\ref{eq2}), (\ref{eq3}) we obtain from (\ref{eq4}) that
$$
A_t=\int_0^{t}\left(Z_s-{\wt U}''(V'(s,X_s)) V'(s,X_s)\lambda_s-\frac{1}{2}{\wt U}'''(V'(s,X_s)) V'(s,X_s)^2\lambda_s\right)^Td\la M\ra_s\lambda_s
$$
\begin{equation}\label{eq5}
- \frac{1}{2}\int_0^{t}\frac{{\wt U}'''(V'(s,X_s)) }{{\wt U}''(V'(s,X_s))^2}d\la N\ra_s.
\end{equation}
Therefore, using  the duality relations
$$
V'(t,X_t)=U'(X_t+Y_t),
$$
$$
 {\wt U}''(V'(t,X_t))=-\frac{1}{U''(X_t+Y_t)},
$$
$$
{\wt U}'''(V'(t,X_t))=-\frac{U'''(X_t+Y_t)}{(U''(X_t+Y_t))^3},
$$
we obtain from   (\ref{eq5})  that
$$
A_t= \int_0^{t}\left( \lambda_s\frac{U'(X_s+Y_s)}{U''(X_s+Y_s)}-
\frac{1}{2}\lambda_s\frac{U'''(X_s+Y_s)U'(X_s+Y_s)^2}{U''(X_s+Y_s)^3}+Z_s\right)^Td\langle M\rangle_s\lambda_s
$$
\begin{equation}\label{eq6}
-\frac{1}{2}\int_0^t  \frac{U'''(X_s+Y_s)}{U''(X_s+Y_s)}d\langle N\rangle_s
\end{equation}
Thus, (\ref{dec}) and (\ref{eq6}) imply that $Y$ satisfies equation (\ref{fb}).

Since
$$
{\wt U}''(V'(s,X_s)) V'(s,X_s)=- \frac{1}{U''(X_s+Y_s)},
$$
 from (\ref{pdex2}) and (\ref{eq2})  we obtain equation (\ref{fbx}) for the optimal wealth.

The proof that  the function $u(t,x)=-{\wt U}'(V'(t,x))- x$
is  the  decoupling field of  the FBSDE (\ref{fb}) is similar. One should take integrals from $s$ to $t$ and use the same arguments.

b)
Since the quadruple  $(Y^{s,x}, Z^{s,x},N^{s,x}, X^{s,x})$ satisfies the FBSDE (\ref{fbs}), (\ref{fbxs}), it follows from the It\^o formula that
for any $t\ge s$
\begin{equation}\label{dc3}
U'(X_t^{s,x}+Y_t^{s,x})=U'(x+u(s,x))- \int_s^{t}\lambda_rU'(X_r^{s,x}+Y_r^{s,x})dM_r
\end{equation}
$$
+\int_s^tU''(X_r^{s,x}+Y_r^{s,x}) dN_r.
$$
Thus  $U'(X_t^{s,x}+Y_t^{s,x}), t\ge s,$ is a local martingale and  a true martingale by assumption. Therefore, it follows from
 (\ref{uth}) and (\ref{utht})  that
\begin{equation}\label{dcu}
 U'(X_t^{s,x}+Y_t^{s,x})=E(U'(X_T^{s,x}+H)/F_t)=V'(t, X_t^{s,x}),
\end{equation}
where the last equality is proved similarly to \cite{S1}.  For $t=s$ we obtain that
\begin{equation}\label{utx2}
U'(x+u(s,x))=V'(s,x),
\end{equation}
 hence
\begin{equation}\label{utx}
u(t,x)=-{\tilde U}'(V'(t,x))-x.
\end{equation}
Since $U(x)$ three-times differentiable and $u(t,x)$ is regular decoupling field,  equality (\ref{utx2}) implies that $V'(t,x)$ will be a regular family of semimartingales.
Therefore, using the It\^o-Ventzel formula for $V'(t, X_t^{s,x})$ and equalities  (\ref{dc3}) , (\ref{dcu}) we have

\begin{equation}\label{dc4}
\int_s^t\big[\varphi'(r, X_r^{s,x})-V''(r, X_r^{s.x})(\lambda_s\frac{U'(X_r^{s,x}+Y_r^{s,x})}{U''(X_r^{s,x}+Y_r^{s,x})}+Z_r^{s,x})\big]dM_r
\end{equation}
$$
+\int_s^tL'(dr,X_r) +\int_s^ta'(r, X_r^{s,x})dK_r
$$
$$
-\int_s^t(\lambda_r\frac{U'(X_r^{s,x}+Y_r^{s,x})}{U''(X_r^{s,x}+Y_r^{s,x})}+Z_r^{s,x})\big)^Td\la M\ra_r(V''(r, X_r^{s,x})\lambda_r+\varphi''(r,X_r^{s,x}))
$$
$$
-\frac{1}{2}\int_s^t(V'''(r,X_r^{s,x}))(\lambda_r\frac{U'(X_r^{s,x}+Y_r^{s,x})}{U''(X_r^{s,x}+Y_r^{s,x})}+Z_r^{s,x}\big)^Td\la M\ra_r(\lambda_r\frac{U'(X_r^{s,x}+Y_r^{s,x})}{U''(X_r^{s,x}+Y_r^{s,x})}+Z_r^{s,x}\big)
$$
$$
= - \int_s^{t}\lambda_rU'(X_r^{s,x}+Y_r^{s,x})dM_r+ \int_s^tU''(X_r^{s,x}+Y_r^{s,x}) dN_r.
$$

Equalizing the integrands of stochastic integrals with respect to $dM$   in (\ref{dc4}) we have that $\mu^{K}$-a.e.
\begin{equation}\label{dc5}
Z_r^{s,x}=\frac{\lambda_rV'(r,X_r^{s,x})+\varphi'(r,X_r^{s,x})}{V''(r,X_r^{s,x})}-\lambda_r\frac{U'(X_r^{s,x}+Y_r^{s,x})}{U''(X_r^{s,x}+Y_r^{s,x})}.
 \end{equation}

Equalizing the parts of finite variations in (\ref{dc4}) taking (\ref{dc5}) in mind we get that for any $t>s$
\begin{equation}\label{dc6}
\int_s^ta'(r, X_r^{s,x})dK_r=\int_s^t\big[\frac{(V''(r, X_r^{s,x})\lambda_r+\varphi''(r,X_r^{s,x}))}{V''(r, X_r^{s,x})}
\end{equation}
$$
-\frac{1}{2}V'''(r,X_r^{s,x})\frac{(V'(r, X_r^{s,x})\lambda_r+\varphi'(r,X_r^{s,x}))}{V''(r, X_r^{s,x})^2}\big]^Td\la M\ra_r(V'(r, X_r^{s,x})\lambda_r+\varphi'(r,X_r^{s,x})).
$$
Let $\tau_s(\varepsilon)=\inf\{t\ge s:K_t-K_s\ge\varepsilon\}$. Since $\la M^i, M^j\ra << \td K$ for any $1\le i,j\le d$, where $\td K=\sum_{i=1}^d\la M^i\ra$, taking
an increasing process $ K+\td K$ (which we denote again by $K$), without loss of generality we can assume that $\la M\ra << K$ and  denote by $C_t$ the matrix of Radon-Nicodym derivatives $C_t=\frac{d\la M\ra_t}{dK_t}$.
Then from (\ref{dc6})

\begin{equation}\label{dc7}
\int_s^{\tau_s(\varepsilon)}\bigg[\frac{(V''(r, X_r^{s,x})\lambda_r+\varphi''(r,X_r^{s,x}))^TC_r(V'(r, X_r^{s,x})\lambda_r+\varphi'(r,X_r^{s,x}))}{V''(r, X_r^{s,x})}
\end{equation}
$$
-\frac{1}{2}V'''(r,X_r^{s,x})\frac{(V'(r, X_r^{s,x})\lambda_r+\varphi'(r,X_r^{s,x}))^TC_r(V'(r, X_r^{s,x})\lambda_r+\varphi'(r,X_r^{s,x}))}{V''(r, X_r^{s,x})^2}
$$
$$
-a'(r, X_r^{s,x})\bigg]dK_r=0.
$$
Since  for any $x\in R$ the process  $X_r^{s,x}$ is continuous function on $\{(r,s), r\ge s\}$ with $X_s^{s,x}=x$ (as a solution of equation (\ref{fbxs})) and $V'(t,x)$ is a regular family of semimartingales,
dividing  equality ( \ref{dc7}) by $\varepsilon$ and passing to the limit as $\varepsilon\to 0$  from \cite{MT3} ( Proposition B1 ) we obtain that for each $x$
\begin{equation}\label{dc8}
a'(s, x)=\frac{(V''(s, x)\lambda_s+\varphi''(s,x))^TC_s(V'(s, x)\lambda_s+\varphi'(s,x))}{V''(s, x)}
\end{equation}
$$
-\frac{1}{2}V'''(s,x)\frac{(V'(s, x)\lambda_s+\varphi'(s,x))^TC_s(V'(s, x)\lambda_s+\varphi'(s,x))}{V''(s, x)^2}
$$
$$
=\frac{1}{2}\bigg(\frac{(V'(s, x)\lambda_s+\varphi'(s,x))^TC_s(V'(s, x)\lambda_s+\varphi'(s,x))}{V''(s,x)}\bigg)',\;\mu^{K} - a.e.,
$$
which implies that $V'(t,x)$ satisfies the BSPDE
$$
    V'(t,x)= V'(0,x) +\frac{1}{2}\int_0^t \bigg(\frac{(V'(s,x)\lambda_s+\varphi'(s,x))^TC_s(V'(s, x)\lambda_s+\varphi'(s,x))} {V''(s,x)}\bigg)'\,dK_s
$$
\begin{equation}\label{pded}
   +\int_0^t\varphi'(s,x)\,dM_s+ L'(t,x), \quad  V'(T,x)=U'(x+H).\qed
\end{equation}

{\bf Remark 1.} In the proof of the part a) of the theorem we need the condition that $V'(t,x)$ is a regular family of semimartingales  only  to show
equality (\ref{op}) and to obtain representation  (\ref{zb}). Equality (\ref{op})  one can prove without this assumption  ( replacing the stochastic line integral by a local martingale
orthogonal to $M$)
from the duality relation
$$
V'(t,X_t(x))=\rho_t(y),\;\;\; y=V'(x),
$$
where $\rho_t(y)/y$ is the density of the minimax martingale measure (see \cite{S3} and \cite{OZ} for the version with random endovment). Since $\rho_t(y)/y$ is representable in the form ${\cal E}(-\lambda\cdot M+D)$, for a local martingale $D$
orthogonal to $M$, using the Dolean Dade equation we have
$$
V'(t,X_t)=\rho_t=y-\int_0^t\lambda_s\rho_sdM_s+\int_0^t\rho_sdD_s=
$$
$$
=1-\int_0^t\lambda_sV'(s,X_s)dM_s+R_t,
$$
where $R_t\equiv (Z\cdot D)_t$ is a local martingale orthogonal to $M$. Further the proof will be the same if we always use a local martingale $R_t$ instead of stochastic line integral $\int_0^t(L'(ds,X_s)$.
Hence the representation  (\ref{zb}) will  be of  the form
$$
N_t=- \int_0^t{\wt U}''({V}'(s,X_s))dR_t.
$$
{\bf Remark 2.} It follows from the proof of Theorem 1, that if a regular decoupling field for the FBSDE (\ref{fb}), (\ref{fbx}) exists, then the
second component of the solution $Z$ is also of the form $Z_t=g(\omega, t, X_t)$ fore some measurable function $g$ and if  we assume that any orthogonal to $M$ local martingale
 $L$ is represented as a stochastic integral with respect to the given
continuous local martingale $M^\bot$, then the third component $N$ of the solution will take the same form $N_t=\int_0^tg^\bot(s,X_s)dM^\bot_s$, for some measurable function $g^\bot$.

\section{Another version of the Forward-Backward system   (\ref{fb})-(\ref{fbx}) }

In this section  we derive other version of the
Forward-Backward system   (\ref{fb}), (\ref{fbx})  and prove an existence of a solution.

{\bf Theorem 2.} Let  utility function $U$ be three-times
continuously differentiable and  let  $S$ be a continuous semimartingale. Assume that conditions (2) and (5) are
satisfied.  Then  there exists a quadruple $(P, \psi, L, X)$ that
satisfies the FBSDE
\begin{equation}\label{nec4}
X_t= x- \int_0^t\frac{\lambda_sP_s+\lambda_sU'(X_s)+\psi_s}{U''(X_s)}dS_s,
 \end{equation}
\be{eqnarray}
\notag P_t=P_0\\
\notag +\int_0^{t}\bigg[\lambda_s
-\frac{1}{2}{U'''(X_s)}\frac{(\lambda_sP_s+\lambda_sU'(X_s)+\psi_s\big)}{U''(X_s)^2}\bigg]^Td\la
M\ra_s\big(\lambda_sP_s+\lambda_sU'(X_s)+\psi_s\big)\\
\lbl{nec3}
+\int_0^t\psi_sdM_s+L_t,\\
\notag P_T=U'(X_T+H)-U'(X_T).
\ee{eqnarray}
In addition the  optimal strategy is expressed as
\begin{equation}\label{nec5}
\pi^*_t=-\frac{\lambda_tP_t+\lambda_tU'(X_t)+\psi_t}{U''(X_t)}
 \end{equation}
and the optimal wealth $X^*$ coincides with $X$.

{\it Proof.}
Define the process
\begin{equation}\label{nec1}
P_t=E(U'( X_T^*+H)/F_t)-{U}'(X^*_t).
\end{equation}
Note that the integrability of  $U'(X_T^*+H)$ follows from the duality relation (\ref{Z}).
It is evident that $P_T=U'(X^*_T+H)-U'(X^*_T)$.

Since $U$ is three-times differentiable, $U'(X^*_t)$ is a continuous semimartingale  and $P_t$ admits the  decomposition
\begin{equation}\label{nec2}
P_t=P_0+A_t+\int_0^t\psi_udM_u+L_t,
\end{equation}
where $A$ is a predictable process of finite variations and $L$ is a local martingale orthogonal to $M$.

Since $\rho^*_t$ is the density of a martingale measure , it is of the form
$\rho^*_t={\cal E}_t(-\lambda\cdot M+R), R\bot M$. Therefore, (\ref{Z}) and (\ref{nec1}) imply that
$$
E(U'( X_T^*+H)/F_t)=y\rho_t^*=y-\int_0^t\lambda_s y\rho^*_sdM_s + \tilde R_t
$$
\begin{equation}\label{nec7}
=y-\int_0^t\big(P_s+U'(X_s^*)\big)\lambda_sdM_s + \tilde R_t,
\end{equation}
 where $y=EU'(X^*_T+H)$  and $\tilde R$ is a local martingale  orthogonal to $M$.

By definition of the process $P_t$, using the It\^o formula for $U'(X_t^*)$ and taking  decompositions (\ref{nec2}), (\ref{nec7})  in mind, we obtain
$$
P_0+A_t+\int_0^t\psi_sdM_s+L_t=y-\int_0^t\big(P_s+U'(X_s^*)\big)\lambda_sdM_s+
\tilde R_t-
$$
$$
-U'(x)
-\int_0^{t}U''(X^*_s)\pi_s^{*T}d\la M\ra_s\ld_s
-\frac{1}{2}\int_0^{t}U'''(X^*_s)\pi_s^{*T}d\la M\ra_s\pi_s^*
$$
 \begin{equation}\label{nec10}
-\int_0^{t}U''(X^*_s)\pi_s^*dM_s.
 \end{equation}
Equalizing the integrands of stochastic integrals with respect to $dM$  we have that $\mu^{\la M\ra}$-a.e.
\begin{equation}\label{nec12}
\pi^*_t=-\frac{\lambda_tP_t+\lambda_tU'(X^*_t)+\psi_t}{U''(X^*_t)}
 \end{equation}
Equalizing the parts of finite variations in (\ref{nec10}) we get
\begin{equation}\label{nec13}
A_t=-\int_0^{t}\big(U''(X^*_s)\ld_s+\frac{1}{2}U'''(X^*_s){\pi_s^*}\big)^Td\la M\ra_s\pi_s^*
\end{equation}
and from (\ref{nec12}), substituting the expression for $\pi^*$ in (\ref{nec13}) we obtain that
\begin{equation}\label{nec14}
A_t=\int_0^{t}\bigg[\lambda_s
-\frac{1}{2}{U'''(X_s)}\frac{(\lambda_sP_s+\lambda_sU'(X_s)+\psi_s\big)}{U''(X_s)^2}\bigg]^Td\la
M\ra_s\big(\lambda_sP_s+\lambda_sU'(X_s)+\psi_s\big)
\end{equation}
Therefore, (\ref{nec14}) and (\ref{nec2}) imply that $P_t$ satisfies equation (\ref{nec3}). Integrating both parts of equality (\ref{nec12}) with respect to $dS$ and adding
the  initial capital we obtain equation (\ref{nec4}) for the optimal wealth.\qed

{\bf Remark 1.} Similarly to Theorem 1b) one can show that  $u(t,x)=V'(t,x)-U'(x)$ is the decoupling field of  (\ref{nec4}),(\ref{nec3}).

{\bf Remark 2.} The generator of equation (\ref{nec3}) does not contain the orthogonal martingale part. Therefore it preserves the same form  without assumption of the continuity of the filtration
(if $S$ is continuous).

{\bf Corollary.} Let conditions of Theorem 2  be satisfied and
assume that the filtration $F$ is continuous.  Then there exists a solution of  FBSDE  (\ref{fb}), (\ref{fbx}).
In particular, if  the pair $(X,P)$
is a solution of  (\ref{nec4}),(\ref{nec3}), then the pair
$(X,Y)$, where
$$
Y_t=- {\tilde U}'(P_t+U'(X_t))-X_t,
$$
satisfies the  FBSDE  (\ref{fb}), (\ref{fbx}).

Conversely, if the pair $(X,Y)$  solves the FBSDE  (\ref{fb}), (\ref{fbx}), then $(X_t,P_t=U'(X_t+Y_t)-U'(X_t))$ satisfies  (\ref{nec4}),(\ref{nec3}).

\end{document}